\documentclass[12]{article} 
\usepackage[T2A]{fontenc}
\usepackage[cp1251]{inputenc}
\usepackage[russian]{babel}
\usepackage{amssymb,amsmath,amsthm,amsfonts}
\usepackage{setspace}
\usepackage{geometry}

\begin{document}
\large
\onehalfspacing
\begin{center}
{\Large {\bf Ре\-гу\-ля\-ри\-зо\-ван\-ный след опе\-ра\-то\-ра Лап\-ла\-са-Бельт\-ра\-ми воз\-му\-щен\-но\-го опе\-ра\-то\-ром умно\-же\-ния на функ\-цию \\на мно\-го\-обра\-зи\-ях со спе\-циаль\-ным воз\-му\-ще\-ни\-ем \\мет\-ри\-ки сфе\-ры}}\\
Т.В.Зыкова\\
16 апреля 2014 г.
\end{center}

Рассмотрим $M$ - двумерное компактное замкнутое многообразие без края такое, что\\
1) $M\in SC_{2\pi}$ \cite{Besse}, т. е. все геодезические $M$ замкнуты и имеют одинаковую длину $2\pi$;\\
2) Метрика $M$ яв\-ляет\-ся воз\-му\-ще\-ни\-ем мет\-ри\-ки стан\-дарт\-ной сфе\-ры $S^2$. Т. е., пусть в $R^3$ есть коор\-ди\-на\-ты $(u_1, u_2, u_3)$ и из\-вест\-на их связь с де\-кар\-то\-вы\-ми коор\-ди\-на\-та\-ми $(x, y, z)$, поз\-во\-ляю\-щая вы\-пи\-сать в этих коор\-ди\-на\-тах ев\-кли\-до\-ву мет\-ри\-ку $ds^2=dx^2+dy^2+dz^2$ в $R^3$, а ог\-ра\-ни\-че\-ни\-ем этой по\-лу\-чив\-шей\-ся мет\-ри\-ки на сфе\-ру еди\-нич\-но\-го ра\-диу\-са (как пра\-ви\-ло за\-ме\-ной од\-ной из пе\-ре\-мен\-ных, на\-при\-мер, $u_3$ на единицу), по\-лу\-чить об\-щий вид стан\-дарт\-ной мет\-ри\-ки дву\-мер\-ной сфе\-ры в коор\-ди\-на\-тах $(u_1, u_2)$:
\begin{equation}
ds^2=A(u_1, u_2) du_1^2+2B(u_1, u_2) du_1 du_2+C(u_1, u_2) du_2^2,
\end{equation}
где  $A(u_1, u_2)$, $B(u_1, u_2)$, $C(u_1, u_2)$ таковы, что $ds^2$ является ри\-ма\-но\-вой и \\$A(u_1, u_2)C(u_1, u_2) -B^2(u_1, u_2)>0$. Тогда под метрикой $M$, будем понимать метрики вида
\begin{equation}
\label{metr}
ds_p^2=A_p(u_1, u_2) du_1^2+2B_p(u_1, u_2) du_1 du_2+C_p(u_1, u_2) du_2^2,
\end{equation}
для которых верно:\\
а) $A_p(u_1, u_2)$, $B_p(u_1, u_2)$, $C_p(u_1, u_2)$ та\-ко\-вы, что $ds^2_p$ яв\-ля\-ет\-ся ри\-ма\-но\-вой, \\$A_p(u_1, u_2)C_p(u_1, u_2) -B_p^2(u_1, u_2)>0$ и гео\-де\-зи\-чес\-кие по\-лу\-чив\-шей\-ся мет\-ри\-ки замк\-ну\-ты и име\-ют оди\-на\-ко\-вую дли\-ну  $2\pi$. \\
б) $A_p(u_1, u_2)=A(u_1, u_2) + P_{A}(u_1, u_2)$, $B_p(u_1, u_2)=B(u_1, u_2) + P_{B}(u_1, u_2)$, \\$C_p(u_1, u_2)=C(u_1, u_2) + P_{C}(u_1, u_2)$, то есть воз\-му\-ще\-ния та\-ко\-вы, что при об\-ну\-ле\-нии функ\-ций $P_{A}(u_1, u_2)$, $P_{B}(u_1, u_2)$, $P_{C}(u_1, u_2)$, мы по\-лу\-чим стан\-дарт\-ную мет\-ри\-ку сфе\-ры $ds^2$. 

Приведем известные примеры таких метрик:\\
(I) Метрика на сфере $ds^2$ заданная в полярных координатах, а $ds_p$ являющаяся мет\-ри\-кой вращения Цолля (см., например, \cite{Besse});\\
(II) Метрика на сфере $ds^2$ заданная в сферо-конических координатах, а $ds_p$ заданная двухпараметрическим семейством $C^\infty$-глад\-ких мет\-рик (попарно не\-изо\-мет\-рич\-ных)  (см. для подробностей \cite{BF1999}).

Мы будем исследовать оператор Лапласа-Бельтрами возмущенный комп\-лекс\-но\-знач\-ной функцией $q$ на многообразиях типа $M$. Следует отметить, что в такой общей пос\-та\-нов\-ке задача нахождения ре\-гу\-ля\-ри\-зо\-ван\-но\-го следа для оператора Лап\-ла\-са-Бельтрами прежде не расс\-мат\-ри\-ва\-лась. В 1991 году В.А.~Садовничим и В.В.~Дуб\-ровс\-ким в \cite{SadDubr1991} для $S^2$ и нечетного потенциала $q$ впервые была получена формула следа со скобками:
$$
\mu_0+\sum\limits_{k=0}^{\infty}\left[\sum\limits_{i=0}^{2k}\mu_{ki}-k(k+1)(2k+1)\right]=-\frac{1}{8\pi}\int\limits_{S^2}q^2dS.
$$ 
Поз\-же В. Е. Подольский \cite{Podolskii1994}, применив к этой задаче суммирование по Абелю и затем к полученной формуле тауберову теорему Литлвуда, доказал, что ряд сходится без скобок (но этот случай является единственным исключением). Поз\-же В.Е. По\-дольс\-ким \cite{P1996} были получены аналогичные формулы для любых сте\-пе\-ней собственных чисел оператора Лапласа-Бельтрами с потенциалом на ком\-пакт\-ных симметрических прост\-ранст\-вах ранга 1. А.Н.~Бобров предпринимал попытку \cite{Bobrov} найти след оператора Лапласа-Бельтрами с потенциалом на по\-верх\-нос\-ти вращения Цолля, но допустил неточность (не была учтена необходимость построения оператора аналогичного $-\tilde{\Delta}_M$ введенного в этой работе, см. ниже) и приведенную им фомулу нельзя считать окончательно верной, но результат для случая простой сферы $S^2$ и произвольной комп\-лекс\-но\-знач\-ной функ\-ции $q\in C^{\infty}$ был получен. В.А.~Са\-дов\-ни\-чий и З.Ю.~Фа\-зул\-лин для оператора возмущенного произвольной комп\-лекс\-но\-знач\-ной функцией лучшей гладкости: в 2005 году для  $q\in C^2(S^2)$\cite{SadFaz2005}, а в 2011 году для $q\in W^2_1(S^2)$\cite{SadFaz2011} получили формулу ре\-гу\-ля\-ри\-зо\-ван\-но\-го следа:
$$
\sum\limits_{k=0}^{\infty}\left[\sum\limits_{i=0}^{2k}\mu_{ki}-k(k+1)(2k+1)-c_0\right]=2c_1,
$$ 
где $c_0=\dfrac{1}{4\pi}\int\limits_{S^2}q(\omega)d\omega$, $c_1=\dfrac{1}{32\pi^3}\int\limits_{S^2}\int\limits_{S^2}\dfrac{q(\omega)q(\omega_0)}{\sqrt{1-(\vec{\omega},\vec{\omega_0})^2}}d\omega d\omega_0 -\dfrac{1}{16\pi}\int\limits_{S^2}q^2(\omega)d\omega$, $\displaystyle (\vec{\omega},\vec{\omega_0})$ - скалярное произведение векторов $\displaystyle \vec{\omega} = (\cos\varphi\sin\theta, \sin\varphi\sin\theta,\cos\theta)$ и \\$\displaystyle \vec{\omega_0} = (\cos\varphi_0\sin\theta_0, \sin\varphi_0\sin\theta_0,\cos\theta_0)$.

В этой работе будет приведена формула следа, об\-об\-щаю\-щая пос\-чи\-тан\-ные ранее частные примеры. 

Кратко обратимся к схеме пост\-рое\-ния ре\-гу\-ля\-ри\-зо\-ван\-но\-го оператора Лапласа-Бельт\-ра\-ми на многообразиях типа $M$.

Для нормализованного оператора Лапласа-Бельтрами на стандартной сфере $S^2$ хо\-ро\-шо известен (см., например, \cite{H19861988}) его спектр, который состоит из точек 
\begin{equation}
\varkappa_{ki} = k(k+1),
\end{equation} 
где $k=0,1,\dots$; $i=1,\dots, N_k$, где кратность $N_k = 2k+1$.  Обозначим через $-\tilde{\Delta}_M$ псевдодифференциальный оператор (далее ПДО) второго порядка на $M$, спектр которого совпадает с $\varkappa_{ki}$.

Пусть $-\Delta_{M}$ - оператор Лапласа-Бельтрами на многообразии $M$. А.~Вейнстейн в \cite{W1977} указал на то, что оператор $-\Delta_{M}$ представим в виде $-\Delta_{M}=-\tilde{\Delta}_M + B$, где $B$ - ПДО нулевого порядка. Такое представление поз\-во\-ли\-ло полагать, что основные результаты этой работы, посвященной исследованию асимп\-то\-тик клас\-те\-ров разностей собственных значений возмущенного и не\-воз\-му\-щен\-но\-го операторов, применимы к собственным чис\-лам опе\-ра\-то\-ра $-\Delta_{M}$. Но для получения окончательного вида асимп\-то\-тик клас\-те\-ров, необ\-хо\-ди\-мо знать сим\-вол оператора $B$ или его усреднение. А.~Вейн\-стейн выдвинул гипотезу, что сим\-волом усреднения сим\-вола такого воз\-му\-ще\-ния яв\-ля\-ет\-ся усреднение гаусовой кривизны многообразия вдоль геодезической. Позже, в 1997 году С. Зельдич \cite{Zelditch1997}, вычислил явный вид сим\-вола усреднения, указав, что помимо пред\-по\-ла\-гае\-мо\-го А. Вейн\-стей\-ном, есть еще один член не всегда равный нулю:
\begin{equation}
\label{meancurvature2}
\dfrac{1}{8\pi}\int\limits_{\gamma}\left(K_M -1 +
\left[\frac{1}{3}(K_M)_\mathrm{v} u^3\int\limits_0^r(K_M)_\mathrm{v}J^3dt - (K_M)_\mathrm{v} u^2 J\int\limits_0^r(K_M)_\mathrm{v}uJ^2dt\right]\right)dr,
\end{equation}
здесь $K_M$ - гаус\-со\-ва кри\-виз\-на мно\-го\-об\-ра\-зия, $\gamma$ --- произ\-воль\-ная гео\-де\-зи\-чес\-кая, $\mathrm{v}$-еди\-нич\-ный век\-тор нор\-ма\-ли к гео\-де\-зи\-чес\-кой $\gamma$, $J(r,\omega)$ - объем\-ная плот\-ность в гео\-де\-зи\-чес\-ких по\-ляр\-ных коор\-ди\-на\-тах, то есть $dvol(\gamma)=J(r,\omega)dr d\omega$, $u$ и $v$ - фун\-да\-мен\-таль\-ные ре\-ше\-ния урав\-не\-ния Якоби вдоль гео\-де\-зи\-чес\-кой $\gamma$ с ус\-ло\-вия\-ми $\left(
\begin{array}{rr}
u(0) & v(0) \\
\dot u(0) & \dot v(0) \\
\end{array}
\right)=
\left(
\begin{array}{rr}
1 & 0 \\
0 & 1 \\
\end{array}
\right)$. \\

Нам удобней рассматривать эту функцию не на пространстве геодезических, а на ко\-ка\-са\-тель\-ном расслоении единичных сфер, поэтому перепишем ее в сле\-дую\-щем виде:
\begin{equation}
\label{meancurvature}
\begin{aligned}
&\sigma^{av}= \dfrac{1}{2\pi}\int\limits_{0}^{2\pi}(\exp t\Xi)^*(\sigma)dt,\quad\text{где}\\
&\sigma = \frac{1}{4}\left(K_M-1 +
\left[\frac{1}{3}(K_M)_\mathrm{v} u^3\int\limits_0^r(K_M)_\mathrm{v}J^3ds - (K_M)_\mathrm{v} u^2 J\int\limits_0^r(K_M)_\mathrm{v}uJ^2ds\right]\right),
\end{aligned}
\end{equation}
и $\Xi$ - гамильтоново векторное поле на ко\-ка\-са\-тель\-ном рас\-слое\-нии $T^*M\setminus \{0\}$.

Обозначим собственные числа $-\Delta_{M}$ через $\lambda_{ki}$, где $k=0,1,\dots$; $i=1,\dots, 2k+1$. Здесь двойная нумерация согласована с нумерацией $\varkappa_{ki}$. Отметим, что для собственных чисел оператора $-\Delta_{M}$ имеет место следующая оценка \cite{W1976}:
\begin{equation}
\label{eigenvalueSpMl}
|\lambda_{ki}-\varkappa_{ki}|=O(1).
\end{equation}

Далее рассмотрим $-\Delta_{M}+q$ ---  оператор  Лапласа-Бельтрами возмущенный комп\-лекс\-но\-значной функцией на $M$ и обозначим собственные числа этого опе\-ра\-то\-ра через $\mu_{ki}$, где двойная нумерация согласована с нумерацией $\lambda_{ki}$ (а значит, и с нумерацией $\varkappa_{ki}$), и также $k=0,1,\dots$; $i=1,\dots, 2k+1$, и также верна аналогичная оценка:
\begin{equation}
\label{eigenvalueMlPertub}
|\mu_{ki} - \lambda_{ki}|=O(1).
\end{equation}

Наша задача получить формулу ре\-гу\-ля\-ри\-зо\-ван\-но\-го следа для собственных чисел оператора $-\Delta_{M}+q$ и $\varkappa_{ki}$  (то есть, нам интересно узнать, как <<отличаются>> собст\-вен\-ные числа оператора Лапласа-Бельтрами с потенциалом взятого на $M$ от собственных чисел оператора Лапласа-Бельтрами рассматриваемого на стан\-дарт\-ной сфере). Ре\-ше\-ние поставленной задачи будет проведено в три этапа:\\
1. Построение формулы ре\-гу\-ля\-ри\-зо\-ван\-но\-го следа для собственных чисел опе\-ра\-то\-ра $-\Delta_{M}=-\tilde{\Delta}_M+B$ и $-\tilde{\Delta}_M$ (то есть для $\lambda_{ki}$ и $\varkappa_{ki}$);\\
2. Построение формулы ре\-гу\-ля\-ри\-зо\-ван\-но\-го следа для собственных чисел опе\-ра\-то\-ра $-\Delta_{M}+q$ и $-\Delta_{M}$ (то есть для $\mu_{ki}$ и $\lambda_{ki}$);\\
3. Сведение и вычисление общей формулы ре\-гу\-ля\-ри\-зо\-ван\-но\-го следа для $\mu_{ki}$ и $\varkappa_{ki}$, с по\-мощью ре\-зуль\-та\-тов пунк\-тов 1. и 2.

Опишем эти этапы более подробно.\\

{\bf 1. Построение формулы ре\-гу\-ля\-ри\-зо\-ван\-но\-го следа для собственных чисел опе\-ра\-то\-ра $-\Delta_{M}=-\tilde{\Delta}_M+B$ и $-\tilde{\Delta}_M$.}
\\

Положим $\nu_{ki}=\lambda_{ki}-k(k+1)$, здесь $i=0,\dots,2k$. Имеет место асимпто\-ти\-чес\-кое разложение при $k \to \infty$ (см., например, \cite{GU1983}):
\begin{equation}
\label{mu}
\sum\limits_{i=0}^{2k}\nu_{ki}= a_0(2k+1)+a_1+a_2(2k+1)^{-1} + O(k^{-2}),
\end{equation}
а значит и
\begin{equation}
\sum\limits_{i=0}^{2k}\lambda_{ki}=(2k+1)k(k+1) + a_0(2k+1)+a_1+a_2(2k+1)^{-1} + O(k^{-2}).
\end{equation}
Из пос\-лед\-ней фор\-му\-лы вид\-но, что для вы\-чис\-ле\-ния пер\-во\-го сле\-да нам необ\-хо\-димо знать коэф\-фи\-циен\-ты $a_0$, $a_1$ и $a_2$.

Введем в рассмотрение ряды $F(t)=\sum\limits_{k=0}^{\infty}(2k+1)e^{-k(k+1) t}$ и $L(t)=\sum\limits_{k=0}^{\infty}e^{-\lambda_k t}$.  Посколь\-ку $\lambda_k\sim k$, то $L(t)$ равномерно сходится при $0<t_0\le t < +\infty$ $\forall t_0$. Исследуя ряд 
$$
L(t) = \sum\limits_{k=0}^{\infty}\sum\limits_{i=0}^{2k}e^{-\lambda_{ki}t} = \sum\limits_{k=0}^{\infty}\sum\limits_{i=0}^{2k}e^{-(k(k+1)+\nu_{ki})t} = \sum_{k=0}^{\infty}e^{-k(k+1)t}\sum\limits_{i=0}^{2k}e^{-\nu_{ki}t},
$$
проведя все необходимые оценки, можно показать, что при $t\to 0$
\begin{equation}
L(t) = F(t)-a_0tF(t)-a_1t\sum_{k=0}^{\infty}e^{-k(k+1)t}-a_2t\sum_{k=0}^{\infty}e^{-k(k+1)t}(2k+1)^{-1} + O(t), 
\end{equation}
а окончательно исследовав асимптотики при $t\to0$ для рядов $\sum_{k=0}^{\infty}e^{-k(k+1)t}$ и \\$\sum_{k=0}^{\infty}e^{-k(k+1)t}(2k+1)^{-1}$  (см., например, \cite{P1996}), получить что 
\begin{equation}
\label{L}
L(t) = F(t)-a_0tF(t) + a_1\frac{1}{2}\sqrt{\pi}t^{1/2} -  a_2t\ln t + O(t).
\end{equation}

Для оператора яв\-ляю\-ще\-го\-ся са\-мо\-соп\-ря\-жен\-ным рас\-ши\-ре\-ни\-ем опе\-ра\-то\-ра $-\Delta_{M}$ из\-вест\-на асимп\-то\-ти\-ка его тета-функ\-ции, сов\-па\-даю\-щая с $L(t)$ ( см., нап\-ри\-мер, \cite{Teilor}):
\begin{equation}
\begin{aligned}
\label{TetL}
&L(t) = t^{-1}\sum\limits_{j=0}^{\infty}l_j t^j,
\end{aligned}
\end{equation}
и аналогично для $F(t)$
\begin{equation}
\begin{aligned}
\label{TetF}
&F(t) =t^{-1}\sum\limits_{j=0}^{\infty}f_j t^j.
\end{aligned}
\end{equation}
Из этих формул и (\ref{L}) следует, что $a_1=0$ и $a_2=0$ и осталось вычислить только $a_0$. 
Подставим (\ref{TetF}) в (\ref{L}) и перепишем:
\begin{equation}
\label{TrGeom2}
L(t) = (f_0\frac{1}{t} + f_1)  - a_0f_0 + O(t),\text{ при }t \to +0.
\end{equation}
Сравнивая получившуюся формулу с (\ref{TetL}), получаем, что\\
$l_0= f_0$ и $l_1 = f_1 - a_0f_0$, а следовательно, $a_0 = \dfrac{f_1-l_1}{f_0}$, где $l_1, f_1, f_0$ коэффициенты соответствующих тета-функций операторов и их можно вычислить через зна\-че\-ния ана\-ли\-ти\-чес\-ко\-го продолжения дзета-функций данных операторов, а значит ряд
$$
\begin{aligned}
\sum\limits_{k=0}^{\infty}\left(\sum_{i=0}^{2k}\lambda_{ki}-k(k+1)(2k+1)-a_0 (2k+1)\right)
\end{aligned}
$$
сходится абсолютно.

Для оценки правой части ряда удобно рассмотреть разность рядов $L(t)-(1-a_0 t)F(t)$ при $t\to0$. Выполнив все необходимы выкладки и оценки можно получить, что
\begin{equation}
\label{endsum}
\begin{aligned}
& \sum\limits_{k=0}^{\infty}\left(\sum\limits_{i=0}^{2k}\lambda_{ki} - k(k+1)(2k+1)-a_0(2k+1)\right)e^{-k(k+1)t} =\\
& =\frac{1}{t}((1-a_0t)F(t)-L(t))+\frac{1}{2}\sum\limits_{k=0}^{\infty}e^{-k(k+1)t}\left(\sum\limits_{i=0}^{2k}(\lambda_{ki}-k(k+1))^2\right)t+o(1).
\end{aligned}
\end{equation}

Рас\-смот\-рим вто\-рое сла\-гае\-мое пра\-вой час\-ти. Как уже отме\-ча\-лось вы\-ше, мо\-жем вос\-поль\-зо\-вать\-ся резуль\-та\-та\-ми ра\-бо\-ты \cite{W1977} для на\-хож\-де\-ния асимп\-то\-ти\-ки клас\-те\-ров, так как сим\-вол ус\-ред\-не\-ния опе\-ра\-то\-ра $B$ нам известен и определен формулой (\ref{meancurvature}). Значит при $k\to\infty$ верно 
$$
\sum\limits_{i=0}^{2k}(\lambda_{ki}-k(k+1))^2\sim \frac{2k+1}{4\pi^2} \int\limits_{S^{*}M}(\sigma^{av})^2dv + O(1),
$$
где $S^{*}M$ расс\-лое\-ние еди\-нич\-ных сфер в ко\-ка\-са\-тель\-ном прост\-ранст\-ве; $dv$ - ка\-но\-ни\-чес\-кий элемент объема $S^{*}M$.

Далее, подставив асимптотики  (\ref{TetL}) и  (\ref{TetF}) в (\ref{endsum}) при $t\to0$ получим

\begin{equation}
\begin{aligned}
\label{otvet}
& \sum\limits_{k=0}^{\infty}\left(\sum_{i=0}^{2k}\lambda_{ki}-k(k+1)(2k+1)-a_0(2k+1)\right)=\\
& = \lim\limits_{t\to +0}\frac{1}{t}((1-a_0t)F(t)-L(t)) + \frac{1}{8\pi^2}\int\limits_{S^{*}M}(\sigma^{av})^2dv =\\
& = f_2 - l_2 - a_0f_1 + \frac{1}{8\pi^2}\int\limits_{S^{*}M}(\sigma^{av})^2dv,
\end{aligned}
\end{equation}
где $\sigma^{av}$ определен в (\ref{meancurvature}), $a_0 = \dfrac{f_1-l_1}{f_0}$, $l_i$ - коэффициенты разложения тета-функ\-ции $L(t)$ (\ref{TetL}) и $f_i$ - коэффициенты разложения тета-функ\-ции $F(t)$ (\ref{TetF}).\\

{\bf 2. Построение формулы ре\-гу\-ля\-ри\-зо\-ван\-но\-го следа для собственных чисел опе\-ра\-то\-ра $-\Delta_{M}+q$ и $-\Delta_{M}$.}\\

Все рассуждения приведенные выше можно повторить и для этого случая -- для пост\-рое\-ния ре\-гу\-ля\-ри\-зо\-ван\-но\-го следа для $\mu_{ki}$ и $\lambda_{ki}$. 
Аналогично положим $\nu'_{ki}=\mu_{ki}-\lambda_{ki}$, здесь $k=1,2,\dots$, $i=0,\dots,2k$ и  запишем асимптотическое разложение при $k \to \infty$:

\begin{equation}
\label{mu2}
\sum\limits_{i=0}^{2k}\nu'_{ki}= b_0(2k+1)+b_1+b_2(2k+1)^{-1} + O(k^{-2}).
\end{equation}

Введем в рассмотрение ряд $M(t)=\sum\limits_{k=0}^{\infty}e^{-\mu_k t}$ и так как для оператора яв\-ляю\-ще\-го\-ся са\-мо\-со\-пря\-жен\-ным рас\-ши\-ре\-ни\-ем опе\-ра\-то\-ра $-\Delta_{M}+q$ из\-вест\-на асимп\-то\-ти\-ка его тета-функ\-ции, сов\-па\-даю\-щая с $M(t)$, то запишем:
\begin{equation}
\begin{aligned}
\label{TetM}
&M(t) = t^{-1}\sum\limits_{j=0}^{\infty}m_j t^j.
\end{aligned}
\end{equation}
Проводя аналогичное исследование, только для собственных чисел $\mu_{ki}$ и $\lambda_{ki}$, с исполь\-зо\-ва\-ни\-ем асимптотик $M(t)$ и $L(t)$ можно показать, что в (\ref{mu2}) $b_1$ и $b_2$  равны нулю, а $b_0 = \dfrac{l_1-m_1}{l_0}$  и ряд 
$$
\begin{aligned}
\sum\limits_{k=0}^{\infty}\left(\sum_{i=0}^{2k}(\mu_{ki}-\lambda_{ki})-b_0 (2k+1)\right)
\end{aligned}
$$
сходится абсолютно.

Проводя аналогичное исследование для нахождения правой части этого ряда, по\-лу\-чим формулу, аналогичную (\ref{endsum}):
\begin{equation}
\label{endsum2}
\begin{aligned}
& \sum\limits_{k=0}^{\infty}\left(\sum\limits_{i=0}^{2k}(\mu_{ki} - \lambda_{ki})-b_0(2k+1)\right)e^{-\lambda_{ki}t} =\\
&=\frac{1}{t}((1-b_0t)L(t)-M(t))+\frac{1}{2}\sum\limits_{k=0}^{\infty}e^{-\lambda_{ki}t}\left(\sum\limits_{i=0}^{2k}(\mu_{ki}-\lambda_{ki})^2\right)t+o(1).
\end{aligned}
\end{equation}

Рассмотрим вто\-рое сла\-гае\-мое пра\-вой час\-ти. Аналогично, из \cite{W1977} сле\-ду\-ет, что для внут\-рен\-не\-го ряда при $k\to\infty$ верно 
$$
\sum\limits_{i=0}^{2k}(\mu_{ki}-\lambda_{ki})^2\sim \frac{2k+1}{4\pi^2} \int\limits_{S^{*}M}(q^{av})^2dv + O(1),
$$
где $S^{*}M$ расс\-лое\-ние еди\-нич\-ных сфер в ко\-ка\-са\-тель\-ном прост\-ранст\-ве; $dv$ - ка\-но\-ни\-чес\-кий элемент объема $S^{*}M$; сим\-вол усреднения
\begin{equation}
\label{qav}
q^{av}=\dfrac{1}{2\pi}\int\limits_0^{2\pi}(\exp t\Xi)^*(q)dt,
\end{equation}
где $\Xi$ - гамильтоново векторное поле на $T^* M\setminus\{0\}$. 

Тогда, окончательно переходя в  (\ref{endsum2}) к пределу при $t\to+0$ и пользуясь фор\-му\-лами (\ref{TetM}) и (\ref{TetL}), получаем:
\begin{equation}
\begin{aligned}
\label{otvet2}
& \sum\limits_{k=0}^{\infty}\left(\sum_{i=0}^{2k}(\mu_{ki}-\lambda_{ki})-b_0(2k+1)\right)=\\
& = \lim\limits_{t\to +0}\frac{1}{t}((1-b_0t)L(t)-M(t)) + \frac{1}{8\pi^2}\int\limits_{S^{*}M}(q^{av})^2dv =\\
& = l_2 -m_2 - b_0l_1 +\frac{1}{8\pi^2}\int\limits_{S^{*}M}(q^{av})^2dv,
\end{aligned}
\end{equation}
где $q^{av}$ определен в (\ref{qav}), $b_0 = \dfrac{l_1-m_1}{l_0}$, $l_i$ - коэффициенты разложения тета-функ\-ции $L(t)$ (\ref{TetL}) и $m_i$ - коэффициенты разложения тета-функ\-ции $M(t)$ (\ref{TetM}).\\
\\
{\bf 3. Сведение и вычисление общей формулы ре\-гу\-ля\-ри\-зо\-ван\-но\-го следа для $-\Delta_{M}+q$}.\\
\\
Сложим левые и правые части формул (\ref{otvet}) и (\ref{otvet2}) и получим:

\begin{equation}
\begin{aligned}
\label{otvet3}
& \sum\limits_{k=0}^{\infty}\left(\sum_{i=0}^{2k}\mu_{ki}-k(k+1)(2k+1)-(a_0+b_0) (2k+1)\right)=\\
& = f_2  - a_0f_1 + \frac{1}{8\pi^2}\int\limits_{S^{*}M}(\sigma^{av})^2dv -m_2 - b_0l_1 +\frac{1}{8\pi^2}\int\limits_{S^{*}M}(q^{av})^2dv,
\end{aligned}
\end{equation}
где $\sigma^{av}$ определен в (\ref{meancurvature}), $q^{av}$ оп\-ре\-де\-лен в (\ref{qav}),  $a_0 = \dfrac{f_1-l_1}{f_0}$, $b_0 = \dfrac{l_1-m_1}{l_0}$,  $m_i$ - коэф\-фи\-циен\-ты раз\-ло\-же\-ния те\-та-функ\-ции $M(t)$ (\ref{TetM}), $l_i$ - коэф\-фи\-циен\-ты раз\-ло\-же\-ния те\-та-функ\-ции $L(t)$ (\ref{TetL}) и $f_i$ - коэф\-фи\-циен\-ты разложения тета-функ\-ции $F(t)$ (\ref{TetF}).

Для получения окончательного ответа нам необходимо предъявить явный вид коэф\-фи\-циен\-тов асимп\-то\-тик $M(t)$, $L(t)$ и $F(t)$, участ\-вую\-щих в ответе. Сна\-ча\-ла от\-ме\-тим, что асимп\-то\-ти\-ку $F(t)$ можно вычислить на\-пря\-мую ис\-сле\-дуя ряд $F(t)=\sum\limits_{k=0}^{\infty}(2k+1)e^{-k(k+1)t}$, например, с помощью \cite{P1996} и получить, что при $t \to +0$
$$
F(t)=t^{-1} + \frac{1}{3} + \frac{1}{15}t+O(t^2),
$$
а значит $f_0=1$, $f_1=\frac{1}{3}$ и $f_2=\frac{1}{15}$. Выше также было отмечено, что  $f_0=l_0=m_0$, а значит $f_0=l_0=m_0=1$.

Для исследования $M(t)$ и $L(t)$ будем пользоваться иными соображениями. Хорошо известно, что $\theta$-функ\-ция и $\zeta$-функ\-ция свя\-за\-ны пре\-об\-ра\-зо\-ва\-нием Мел\-ли\-на и искомые коэффициенты разложения тета-функций  выражаются через аналитические зна\-че\-ния дзета-функций, а именно
$$
m_1 = \zeta_{-\Delta_{M}+q}(0),\quad m_2 = - \zeta_{-\Delta_{M}+q}(1),
$$
где $\zeta_{-\Delta_{M}+q}(0)$ и $\zeta_{-\Delta_{M}+q}(1)$ - зна\-че\-ния аналитического продолжения дзета-функ\-ции опе\-ра\-то\-ра $-\Delta_{M}+q$ в точках $0$ и $1$. Аналогично для $L(t)$:  
$$
\label{lifi}
l_1 = \zeta_{-\Delta_{M}}(0), \quad l_2 = - \zeta_{-\Delta_{M}}(1),
$$
где $\zeta_{-\Delta_{M}}(0)$ и $\zeta_{-\Delta_{M}}(1)$ - зна\-че\-ния аналитического продолжения дзета-функ\-ции опе\-ра\-то\-ра $-\Delta_{M}$ в точках $0$ и $1$.

В теории ПДО известна стандартная процедура (см., например, \cite{Sh1978}) пост\-рое\-ния параметрикса для  классического эллиптического ПДО на замкнутом многообразии. В ее основе --- нахождение сим\-вола параметрикса, как асимп\-то\-ти\-чес\-кой суммы неких однородных функций, которые, в свою очередь опре\-де\-ляют\-ся через од\-но\-род\-ные ком\-по\-нен\-ты сим\-вола самого оператора с при\-ме\-не\-ни\-ем формулы композиции. В случае, когда оператор дифференциальный, с помощью такого метода, удается по\-лу\-чить компоненты разложения сим\-вола резольвенты по рекуррентным фор\-му\-лам и с по\-мощью их, по\-лу\-чить искомые зна\-че\-ния ана\-ли\-ти\-чес\-ко\-го продолжения $\zeta$-функ\-ции. В на\-шем случае, когда мет\-ри\-ка многообразия задана в абстрактном виде, вы\-чис\-ли\-тель\-ная слож\-ность ре\-ше\-ния задачи сильно возросла и ре\-ше\-ние во многом бы\-ло по\-лу\-че\-но с помощью сим\-воль\-но\-го прог\-рам\-ми\-ро\-ва\-ния в пакете Wolfram Mathematica 9 \cite{math}.
$$
\begin{aligned}
&\zeta_{-\Delta_{M}+q}(0)=\frac{1}{4\pi}\iint\limits_{M}\left(\frac{1}{3}K_{M}-q(u_1,u_2)\right)\sqrt{\text{det}g}du_1du_2,\\
\label{zeta01}
&\zeta_{-\Delta_{M}+q}(1)=-\frac{1}{60\pi}\iint\limits_{M}(\Delta_{M} K_{M}+K_{M}^2)\sqrt{\text{det}g}du_1du_2-\\ 
&-\frac{1}{24\pi}\iint\limits_{M}\left(-\Delta_{M} q(u_1,u_2) + 3q^2(u_1,u_2)-2q(u_1,u_2)\gamma(ML)\right)\sqrt{\text{det}g}du_1du_2,
\end{aligned}
$$
где $\sqrt{\text{det}g} =\sqrt{A_p(u_1,u_2)C_p(u_1,u_2)-B_p^2(u_1,u_2)}$ - ко\-рень из оп\-ре\-де\-ли\-те\-ля мат\-ри\-цы мет\-ри\-чес\-ко\-го тен\-зо\-ра и

\begin{spacing}{2.0}
\noindent\({K_M = \dfrac{1}{4 \left(B_p(u_1,u_2)^2-A_p(u_1,u_2) C_p(u_1,u_2)\right)^2}\times}\\
{\left(C_p(u_1,u_2) {A_p}^{'}_{u_2}(u_1,u_2)^2-2 B_p(u_1,u_2) {A_p}^{'}_{u_2}(u_1,u_2) {B_p}^{'}_{u_2}(u_1,u_2)+\right.}\\
{ {A_p}(u_1,u_2)  {A_p}^{'}_{u_2}(u_1,u_2) {C_p}^{'}_{u_2}(u_1,u_2)+2 B_p(u_1,u_2)^2 {A_p}^{''}_{u_2u_2}(u_1,u_2)-}\\
{2  {A_p}(u_1,u_2) {C_p}(u_1,u_2)  {A_p}^{''}_{u_2u_2}(u_1,u_2)-2 {C_p}(u_1,u_2) {B_p}^{'}_{u_2}(u_1,u_2)
 {A_p}^{'}_{u_1}(u_1,u_2)+}\\
{ {B_p}(u_1,u_2) {C_p}^{'}_{u_2}(u_1,u_2)  {A_p}^{'}_{u_1}(u_1,u_2)+4  {B_p}(u_1,u_2) {B_p}^{'}_{u_2}(u_1,u_2)
 {B_p}^{'}_{u_1}(u_1,u_2)-}\\
{2  {A_p}(u_1,u_2) {C_p}^{'}_{u_2}(u_1,u_2)  {B_p}^{'}_{u_1}(u_1,u_2)- {B_p}(u_1,u_2)  {A_p}^{'}_{u_2}(u_1,u_2){C_p}^{'}_{u_1}(u_1,u_2)+}\\
{{C_p}(u_1,u_2) {A_p}^{'}_{u_1}(u_1,u_2) {C_p}^{'}_{u_1}(u_1,u_2)-2  {B_p}(u_1,u_2)  {B_p}^{'}_{u_1}(u_1,u_2)
C^{'}_{u_1}(u_1,u_2)+}\\
{ {A_p}(u_1,u_2){C_p}^{'}_{u_1}(u_1,u_2)^2-4  {B_p}(u_1,u_2)^2  {B_p}^{''}_{u_1u_2}(u_1,u_2)+}\\
{4 A {A_p}(u_1,u_2) {C_p}(u_1,u_2)  {B_p}^{''}_{u_1u_2}(u_1,u_2)+2  {B_p}(u_1,u_2)^2 {C_p}^{''}_{u_1u_1}(u_1,u_2)-}\\
{\left.2  {A_p}(u_1,u_2){C_p}(u_1,u_2){C_p}^{''}_{u_1u_1}(u_1,u_2)\right)\text{ --- гаус\-со\-ва кри\-виз\-на }M.}\)
\end{spacing}

Таким образом, получили, что 
\begin{equation}
\begin{aligned}
&m_1=\frac{1}{4\pi}\iint\limits_{M}\left(\frac{1}{3}K_{M}-q(u_1,u_2)\right)\sqrt{\text{det}g}du_1du_2,\\
&m_2=\frac{1}{60\pi}\iint\limits_{M}(\Delta_{M} K_{M}+K_{M}^2)\sqrt{\text{det}g}du_1du_2+\\ 
&+\frac{1}{24\pi}\iint\limits_{M}\left(-\Delta_{M} q(u_1,u_2) + 3q^2(u_1,u_2)-2q(u_1,u_2)\gamma(ML)\right)\sqrt{\text{det}g}du_1du_2,
\end{aligned}
\end{equation}
Заметим, что вычисление этих значений опиралось в большей степени на вид метрики многообразия и на вид оператора на ней. Нетрудно увидеть, что  \\$\zeta_{-\Delta_{M}+q}(*)$ при $q=0$, есть $\zeta_{-\Delta_{M}}(*)$. Тогда получаем, что
\begin{equation}
\begin{aligned}
&l_1=\frac{1}{4\pi}\iint\limits_{M}\left(\frac{1}{3}K_{M}\right)\sqrt{\text{det}g}du_1du_2, 
&l_2=\frac{1}{60\pi}\iint\limits_{M}(\Delta_{M} K_{M}+K_{M}^2)\sqrt{\text{det}g}du_1du_2.\\ 
\end{aligned}
\end{equation}
Формулы содержащие интеграл от гауссовой кривизны по многообразию (для $m_1$ и $l_1$) можно упростить используя формулу Гаусса-Бонне, то есть воспользоваться тем, что $\iint\limits_{M}K_{M}\sqrt{\text{det}g}dv_2dv_3= 2\pi \chi (M)$, где $\chi(M)$ - характеристика Эйлера, которая в данном случае равна $2$. И тогда окончательно можем сформулировать главный результат этой работы.\\
\\
{\bf ТЕОРЕМА:} Пусть $M\in SC_{2\pi}$ --- многообразие, метрика которого является воз\-му\-ще\-ни\-ем метрики стандартной сферы и задана формулой (\ref{metr}), $q$ - бесконечно-диф\-фе\-рен\-ци\-руе\-мая комп\-лекс\-но\-знач\-ная функция на $M$, тогда для собственных чисел оператора $-\Delta_{M}+q$ верно равенство:
$$
 \sum\limits_{k=0}^{\infty}\sum\limits_{i=0}^{2k}\left(\mu_{ki}-k(k+1)-\frac{1}{4\pi}\int\limits_{M}qdS\right)=
$$
$$
=\frac{1}{8\pi^2}\int\limits_{S^{*}M}(q^{av})^2dv+\frac{1}{8\pi^2}\int\limits_{S^{*}M}(\sigma^{av})^2dv+\dfrac{1}{15}-
$$
$$
- \frac{1}{60\pi}\int\limits_{M}(\Delta_{M} K_{M}+K_{M}^2)dS-\frac{1}{24\pi}\int\limits_{M}\left(-\Delta_{M} q + 3q^2-2q(K_{M}-1)\right)dS,
$$
где $K_M$ - гаус\-со\-ва кри\-виз\-на $M$, $S^{*}M$ - рас\-слое\-ние еди\-нич\-ных косфер над $M$, $dv$ - ка\-но\-ни\-чес\-кая фор\-ма объема на $S^{*}M$, $q^{av}=\dfrac{1}{2\pi}\int\limits_0^{2\pi}(\exp t\Xi)^*(q)dt$, где $\Xi$ - га\-миль\-то\-но\-во векторное поле на ко\-ка\-са\-тель\-ном рас\-слое\-нии $T^*M\setminus \{0\}$, оп\-ре\-де\-ляе\-мое ри\-ма\-но\-вой струк\-ту\-рой на $M$, $\sigma^{av}= \dfrac{1}{2\pi}\int\limits_{0}^{2\pi}(\exp t\Xi)^*(\sigma)dt$, где $\sigma = \dfrac{1}{4}(K_M-1 +$\\
$\left[\frac{1}{3}(K_M)_\mathrm{v} u^3\int\limits_0^r(K_M)_\mathrm{v}J^3ds - (K_M)_\mathrm{v} u^2 J\int\limits_0^r(K_M)_\mathrm{v}uJ^2ds\right])$, где $\mathrm{v}$ - еди\-нич\-ный вектор нор\-ма\-ли к гео\-де\-зи\-чес\-кой $\gamma$, $J(r,\omega)$ - объемная плотность в гео\-де\-зи\-чес\-ких по\-ляр\-ных коор\-ди\-на\-тах, то есть $dvol(\gamma)=J(r,\omega)dr d\omega$, $u$ и $v$ - фун\-да\-мен\-таль\-ные решения уравнения Якоби вдоль гео\-де\-зи\-чес\-кой $\gamma$ с ус\-ло\-вия\-ми $\left(
\begin{array}{rr}
u(0) & v(0) \\
\dot u(0) & \dot v(0) \\
\end{array}
\right)=
\left(
\begin{array}{rr}
1 & 0 \\
0 & 1 \\
\end{array}
\right)$. 
\newpage{}

\noindent
\begin{tabular}{l p{100pt} r}
 & &\\
~~Т.В. Зыкова & &\\
Механико-математический факультет & &\\
Московского государственного университета имени М.В.Ломоносова, & &\\
zytanya@yandex.ru & &\\
\end{tabular}

\end{document}